\documentclass[a4paper,11pt]{amsart}

\usepackage{amssymb,amsmath,epsfig,graphicx,eufrak}
\usepackage{subfigure}
\usepackage{epsfig,stmaryrd}

\author[Florent Benaych-Georges]{Florent Benaych-Georges}\address{Florent Benaych-Georges, LPMA,  UPMC Univ Paris 6, Case courier 188, 4, Place Jussieu, 75252 Paris Cedex 05, France\\ and CMAP, \'Ecole Polytechnique, route de Saclay, 91128 Palaiseau Cedex, France} \email{florent.benaych@gmail.com}
\title[On a surprising relation between the rectangular and square free convolutions]{On a surprising relation between the Marchenko-Pastur law,  rectangular and square free convolutions}
\date{\today}
\newcommand{\MP}{\operatorname{MP}}

\newcommand{\be}{\begin{equation}}
\newcommand{\ee}{\end{equation}}
\newcommand{\beq}{\begin{eqnarray*}}
\newcommand{\eeq}{\end{eqnarray*}}

\newcommand{\NC}{\operatorname{NC}}

\newcommand{\ssi}{if and only if }

\newcommand{\R}{\mathbb{R}}
\newcommand{\C}{\mathbb{C}}

\newcommand{\ud}{\mathrm{d}}

\newcommand{\pro}{probability }

\newcommand{\f}{\frac}
\newcommand{\ff}{\frac{1}}
\newcommand{\lf}{\left}
\newcommand{\ri}{\right}

\newcommand{\st}{such that }
\newcommand{\la}{\lambda}

\newcommand{\ste}{\, ;\, }

\newcommand{\arc}{\boxplus_{\la}}

\newcommand{\bxp}{\boxplus}
\newcommand{\bxt}{\boxtimes}

\newcommand{\bck}{\backslash}

\newtheorem{Th}{Theorem}
\newtheorem{propo}[Th]{Proposition} 
 
\newtheorem{lem}[Th]{Lemma}

\newtheorem{rmq}[Th]{Remark}
\newtheorem{cor}[Th]{Corollary}

\newenvironment{pr}{\noindent {\bf Proof. }}{\ \ \  \hfill $\square$}

\long\def\symbolfootnote[#1]#2{\begingroup
\def\thefootnote{\fnsymbol{footnote}}\footnote[#1]{#2}\endgroup}


 %Random matrices 

\thanks{This work was partially supported by the \emph{Agence Nationale de la
Recherche} grant ANR-08-BLAN-0311-03.}
\begin{document}
\maketitle

\begin{abstract}In this paper, we prove a result linking the square and the rectangular $R$-transforms, the consequence of which is a surprising relation between the square and rectangular versions  the free additive  convolutions, involving the Marchenko-Pastur law. Consequences on random matrices, on infinite divisibility and on the arithmetics of the square versions of the free additive and multiplicative convolutions  are given.
\\
\\
Dans cet article, on prouve un r\'esultat reliant les versions carr\'e et rectangulaire de la $R$-transform\'ee, qui a pour cons\'equence une relation surprenante entre les versions carr\'e et rectangulaire  de la convolution libre additive, impliquant la loi de Marchenko-Pastur. On donne des cons\'equences de ce r\'esultat portant sur les matrices al\'eatoires, sur l'infinie divisibilit\'e et sur l'arithm\'etique des versions carr\'e des convolutions additives et multiplicatives.
\end{abstract}
%arXiv abtract: In this paper, we prove a result linking the square and the rectangular R-transforms, the consequence of which is a surprising relation between the square and rectangular free convolutions, involving the Marchenko-Pastur law. Consequences on random matrices, on infinite divisibility and on the arithmetics of Voiculescu's free additive and multiplicative convolutions  are given.

\tableofcontents

\section*{Introduction}Free convolutions are operations on \pro measures on the real line which allow to compute the empirical spectral\footnote{The {\em empirical spectral measure} of a matrix is the uniform law on its eigenvalues with multiplicity.}  or singular\footnote{The {\em empirical singular measure} of a matrix $M$ with size $n$ by $p$ ($n\leq p$) is the 
empirical spectral measure of $|M|:=\sqrt{MM^*}$.}  measures of large random matrices which are expressed as sums or products of independent random matrices, the spectral measures of which are known. 
More specifically, the operations $\bxp,\bxt$, called respectively {\em free additive and multiplicative convolutions} are defined in the following way \cite{vdn91}. Let, for each $n$, $M_n$, $N_n$ be $n$ by $n$ independent random hermitian matrices, one of them having a distribution which is invariant under the action of the unitary group by conjugation, the empirical spectral measures of which converge, as $n$ tends to infinity, to non random \pro measures denoted respectively  by $\tau_1, \tau_2$. Then $\tau_1\bxp\tau_2$ is the limit of the empirical spectral law of $M_n+N_n$ and, in the case where the matrices are positive, $\tau_1\bxt\tau_2$ is the limit of the empirical spectral law of $M_nN_n$. In the same way, for any $\la\in [0,1]$, the {\em rectangular free convolution} $\arc$ is defined, in  \cite{bg07}, in the following way. Let $M_{n,p}, N_{n,p}$ be $n$ by $p$ independent random  matrices, one of them having a distribution which is   invariant by multiplication by any unitary matrix on any side,    the symmetrized\footnote{The {\em symmetrization} of a law $\mu$ on $[0,+\infty)$ is the law $\nu$ defined by $\nu(A)=\f{\mu(A)+\mu(-A)}{2}$ for all Borel set $A$.   Dealing with laws on $[0,+\infty)$ or with their symmetrizations is equivalent, but for historical reasons, the rectangular free convolutions have been defined with symmetric laws. In all this paper, we shall often pass from symmetric laws to laws on $[0,+\infty)$ and vice-versa. Thus in order to avoid confusion, we shall mainly use the letter $\mu$ for laws on $[0,\infty)$ and $\nu$ for symmetric ones.} empirical singular measures of which tend, as $n,p$ tend to infinity in such a way that $n/p$ tends to $\la$,  to non random \pro measures $\nu_1,\nu_2$.  Then the symmetrized  empirical singular law of $M_{n,p}+N_{n,p}$ tends to  $\nu_1\arc \nu_2$. These operations can also, equivalently, be defined   in reference to free elements of a non commutative \pro space, but in this paper, we have chosen to use the random matrix point of view.  

 In the cases $\la=0$ or  $\la=1$, i.e. where the rectangular random matrices considered in the previous definition are either ``almost flat" or ``almost square",  the rectangular free convolution with ratio $\la$ can be expressed with the additive free convolution:  $\bxp_1=\bxp$ and  for all symmetric laws $\nu_1,\nu_2$, $\nu_1\bxp_0 \nu_2$ is the symmetric law the push-forward by    the map $t\mapsto t^2$ of which is the free additive convolution of the push forwards of $\nu_1$ and $\nu_2$ by the same map. These surprising relations have no simple explanations, but they allow to hope a general relation between the operations $\arc$ and $\bxp$, which would be true for any $\la$. Up to now, despite many efforts, no such relation had been found, until a paper of Debbah and Ryan \cite{dr07}, where a relation between $\arc,\bxp$ and $\bxt$ is proved in a particular case. In the present paper, we give a shorter proof of a wide generalization\footnote{See Remark \ref{28.08.08.1}.} of their result:  for any $\la\in (0,1]$, we define $\mu_\la$ to be the law of $\la$ times a random variable with law the Marchenko-Pastur 
law with mean $1/\la$, and we prove that for any pair $\mu,\mu'$ of \pro measures on $[0,+\infty)$, we have 
\be\label{27.08.08.77}\sqrt{\mu\bxt\mu_\la}\arc \sqrt{\mu'\bxt\mu_\la}=\sqrt{(\mu\bxp\mu')\bxt\mu_\la}, \ee where for any \pro measure $\rho$ on $[0,+\infty)$, $\sqrt{\rho}$ denotes the symmetrization of the push-forward of $\rho$ by the map $t\mapsto \sqrt{t}$.  
Our proof is based on the following relation between the $R$-transform\footnote{Note that there are two conventions regarding the $R$-transform. The one we use is the one used in the combinatorial approach to freeness  \cite{ns06}, which is not exactly the one used in the analytic approach \cite{hiai}: $R_\mu^{\textrm{combinatorics}}(z)=zR_\mu^{\textrm{analysis}}(z)$.\label{footnote.25.08.08.1}}  $R_\mu$ of a \pro measure $\mu$ on $[0,+\infty)$ and the rectangular $R$-transform $C_{\sqrt{\mu\bxt\mu_\la}}$ with ratio $\la$ of $\sqrt{\mu\bxt\mu_\la}$: we prove that for all $z$, $$R_\mu(z)=C_{\sqrt{\mu\bxt\mu_\la}}(z).$$This relation also allows us to   prove precise relations between $\bxp$-infinitely divisible laws and $\arc$-infinitely divisible laws. 

%We would like to observe that formula \eqref{27.08.08.77} has some consequences which are far from obvious. It means that for $n,p$ large integers \st $n/p\simeq \la$, for $A,B,M$ independent random matrices with respective sizes $n\times n$, $n\times n$ and $p\times p$ \st$A,B$ are invariant in law under left and right multiplication by unitary matrices and $M$ has independent Gaussian entries, if ones defines $P$ to be the $n\times p$ matrix $P=\begin{bmatrix}I_n &0\end{bmatrix}$, then 	as far as the spectral measure is concerned, $$APM(APM)^*+BPM(BPM)^*\simeq (APM+BPM)(APM+BPM)^*.$$It also means, if $1<< n<< p$, that for $M,N$ independent $n\times p$ random matrices, as far as the spectrums are concerned, $$ MM^*+NN^*\simeq (M+N)(M+N)^*.$$

We would like to observe that formula \eqref{27.08.08.77} has some consequences which are far from obvious. It means that for $n,p$ large integers \st $n/p\simeq \la$, for $A,B,M,M'$ independent random matrices with respective sizes $n\times n$, $n\times n$, $n\times p$ and $n\times p$ \st $A,B$ are invariant in law under left and right multiplication by unitary matrices and $M,M'$ have independent Gaussian entries,   then 	as far as the spectral measure is concerned, $$ (AM+BM')(AM+BM')^*\simeq AM(AM)^*+BM'(BM')^*.$$ It also means, if $1<< n<< p$, that for $X,Y$ independent $n\times p$ random matrices, as far as the spectrums are concerned, $$ (X+Y)(X+Y)^*\simeq XX^*+YY^*.$$

The relation \eqref{27.08.08.77} has also  consequences on the arithmetics of   free additive and multiplicative convolutions $\bxp$ and $\bxt$ (Corollaries  \ref{27.08.08.1cor} and \ref{28.08.08.randocycle})
   which wasn't known yet, despite the many papers written the last years about questions related to this subject, e.g.  \cite{bv95,appenice,fbg04,gotzea, gotzeb,bbg07,bbcc08}.
   
{\bf Acknowledgments:} The author would like to thank his friend Raj Rao for bringing the paper \cite{dr07} to his attention and %M\'erouane Debbah,  
{\O}yvind Ryan and Serban Belinschi for some useful discussions. He would also like to thank an anonymous referee for suggesting him Remark \ref{10.6.09.1}. 

\section{A relation between the Marchenko-Pastur law, the square and the rectangular free convolutions}
\subsection{Prerequisites on square and rectangular analytic transforms of \pro measures}
\subsubsection{The square case: the $R$- and $S$-transforms}
 These are analytic transforms of \pro measures which allow to compute the operations $\bxp$ and $\bxt$, like the Fourier transform for the classical convolution. The $R$-transform can be defined for any \pro measure on the real line, but we shall only define it for \pro measures on $[0,+\infty)$. Consider such a \pro measure $\mu$. It $\mu=\delta_0$, then $R_\mu=S_\mu=0$. Now, let us suppose that $\mu\neq \delta_0$. Let us define the function $$M_\mu(z)=\int_{t\in\R}\f{tz}{1-tz}\ud \mu(t).$$ Then the {\it $R$- and $S$-transforms\footnote{See the footnote \ref{footnote.25.08.08.1}.}} of $\mu$, denoted respectively by $R_\mu$ and $S_\mu$ are the analytic functions defined as follows\be\label{23.08.08.1}R_\mu(z)=[(1+z)M_\mu^{-1}(z)]^{-1},\quad S_\mu(z)=\f{1+z}{z}M_\mu^{-1}(z),\ee where the exponent $^{-1}$ refers to the inversion of functions with respect to the operation of composition $\circ$. Note that $M_\mu$ is an analytic function defined in $\{z\in \C\ste 1/z\notin \operatorname{support}(\mu)\}$. Hence in the case where $\mu$ is compactly supported, the functions $M_\mu$ and $(1+z)M_\mu^{-1}(z)$ can be inverted in a neighborhood of zero as  
analytic  functions in a neighborhood of zero vanishing at zero, with non null derivative at zero. In the case where $\mu$ is not compactly supported, these functions are inverted as functions on intervals $(-\epsilon, 0)$ which are equivalent to $(\textrm{positive constant})\times z$ at zero \cite{defconv}.

Note that puting together both equations of \eqref{23.08.08.1}, one gets \be\label{23.08.08.2}S_\mu(z)=\ff{z}R_\mu^{-1}(z)=\f{1+z}{z}M_\mu^{-1}(z).
\ee

The main properties of the $R$- and $S$-transforms are the fact that they characterize measures and their weak convergence and that they allow to compute free convolutions : for all $\mu,\nu$, \be\label{25.08.08.1}R_{\mu\bxp\nu}=R_\mu+R_\nu\quad\textrm{ and }\quad S_{\mu\bxt\nu}=S_\mu S_\nu.\ee

\subsubsection{The rectangular case: the rectangular $R$-transform with ratio $\la$}In the same way, for $\la\in [0,1]$, the rectangular free convolution with ratio $\la$ can be computed with an analytic transform of \pro measures. Let $\nu$ be a symmetric \pro measure on the real line. Let us define 
$H_\nu(z)= z(\la M_{\nu^2}(z)+1)(M_{\nu^2}(z)+1)$, where $\nu^2$ denotes the push forward of $\nu$ by the map $t\mapsto t^2$. Then with the same conventions about inverses of functions than in the previous section,  the {\em rectangular $R$-transform with ratio $\la$} of $\nu$ is defined to be $$C_\nu(z)=U\lf( \f{z}{H_\nu^{-1}(z)}-1\ri), $$where $U(z)=  \f{-\la-1+\lf[(\la+1)^2+4\la z\ri]^{1/2}}{2\la}$ for $\la>0$ and $U(z)=z$ for $\la=0$. By Theorems 3.8, 3.11 and 3.12 of  \cite{bg07}, the  rectangular $R$-transform characterizes measures and their weak convergence, and for all pair $\nu_1, \nu_2$ of symmetric \pro measures, $\nu_1\arc\nu_2$ is characterized by the fact that \be\label{28.08.08.17}C_{\nu_1\arc\nu_2}=C_{\nu_1}+C_{\nu_2}.\ee

\subsection{A relation between the square and the rectangular $R$-transforms}Let us fix $\la\in [0,1]$. We recall that for any \pro measure $\rho$ on $[0,+\infty)$, $\sqrt{\rho}$ denotes the symmetrization of the push-forward of $\rho$ by the map $t\mapsto \sqrt{t}$ and that  for $\la>0$, we have defined $\mu_\la$ to be the law of $\la$ times a random variable with law the Marchenko-Pastur 
law with mean $1/\la$, i.e. $\mu_\la$ is the law with support $[(1-\sqrt{\la})^2, (1+\sqrt{\la})^2]$ and density $$x\mapsto \f{\sqrt{4\la -(x-1-\la)^2}}{2\pi\la x}.$$
For $\la=0$, we let $\mu_\la$ denote the Dirac mass at $1$. 

\begin{Th}\label{R=C}Let $\mu$ be a \pro measure on $[0,+\infty)$. Then we have $$R_\mu(z)=C_{\sqrt{\mu\bxt\mu_\la}}(z).$$
\end{Th}

\begin{rmq}[The cumulants point of view]\label{28.08.08.16} {\rm Suppose  $\mu$   to be compactly supported. Let us denote  the free cumulants \cite{ns06} of $\mu$ by $(k_n(\mu))_{n\geq 1}$ and the rectangular free cumulants with ratio $\la$ \cite{bg07, bg07c} of $\sqrt{\mu\bxt\mu_\la}$ by $(c_{2n}(\sqrt{\mu\bxt\mu_\la}))_{n\geq 1}$. Then the previous theorem means that for all $n\geq 1$, one has \be\label{10.6.09.2}k_n(\mu)=c_{2n}(\sqrt{\mu\bxt\mu_\la}).\ee}
\end{rmq}

\begin{pr} - First of all, note that by continuity of the applications $\mu\mapsto\mu\bxt\mu_\la$, $\rho\mapsto R_\rho$ and $\nu\mapsto C_\nu$ with respect to weak convergence \cite{defconv, bg07}, it suffices to prove the result in the case where $\mu$ is compactly supported. In this case, the functions $M_{\mu}, R_\mu, S_\mu, M_{\mu\bxt\mu_\la}, H_{{\sqrt{\mu\bxt\mu_\la}}},C_{\sqrt{\mu\bxt\mu_\la}}$  are analytic in a neighborhood of zero and the operations of inversion on these functions or related ones can be used without precaution. 

- If $\la>0$, the free cumulants of the Marchenko-Pastur law with mean $1/\la$ are all equal to $1/\la$, thus  the ones of $\mu_\la$ are given by the formula $k_n(\mu_\la)=\la^{n-1}$ for all $n\geq 1$ and $R_{\mu_\la}(z)=\sum_{n\geq 1} \la^{n-1}z^n$. From \eqref{23.08.08.2}, it follows that $S_{\mu_\la}(z)=\ff{1+\la z}$. Hence by \eqref{25.08.08.1}, we have $S_{\mu\bxt\mu_\la}(z)=\f{S_\mu(z)}{1+\la z}$, and by \eqref{23.08.08.2}, \be\label{26.08.08.1}M_{\mu\bxt\mu_\la}(z)=\lf(\f{M_\mu^{-1}(z)}{1+\la z}\ri)^{-1}.\ee Note that since $\mu_0=\delta_1$, \eqref{26.08.08.1} stays true if $\la=0$. Now, let us define the function $T(x)=(\la x+1)(x+1)$. Note that $T(U(x-1))=x$ for $x$ in a neighborhood of zero.
We have $$H_{\sqrt{\mu\bxt\mu_\la}}(z)=z\times T\circ M_{\mu\bxt\mu_\la}(z)=z\times T\circ \lf(\f{M_\mu^{-1}(z)}{1+\la z}\ri)^{-1},$$ and \be\label{28.08.08.4}C_{\sqrt{\mu\bxt\mu_\la}}(z)=U\lf(\f{z}{\lf(z\times T\circ \lf(\f{M_\mu^{-1}(z)}{1+\la z}\ri)^{-1}\ri)^{-1}}-1\ri).\ee

- Hence by \eqref{23.08.08.1} and \eqref{28.08.08.4},  we have the following equivalence\begin{eqnarray*}R_\mu=C_{\sqrt{\mu\bxt\mu_\la}}&\iff& \lf((z+1)M_\mu^{-1}(z)\ri)^{-1}=U\lf(\f{z}{\lf(z\times T\circ \lf(\f{M_\mu^{-1}(z)}{1+\la z}\ri)^{-1}\ri)^{-1}}-1\ri)\\ &\iff&
T\circ  \lf((z+1)M_\mu^{-1}(z)\ri)^{-1}=\f{z}{\lf(z\times T\circ \lf(\f{M_\mu^{-1}(z)}{1+\la z}\ri)^{-1}\ri)^{-1}}\\
&\iff& \lf(z\times T\circ \lf(\f{M_\mu^{-1}(z)}{1+\la z}\ri)^{-1}\ri)^{-1}\times  T\circ  \lf((z+1)M_\mu^{-1}(z)\ri)^{-1}=z. 
\end{eqnarray*}Composing both terms on the right by $(z+1)M_\mu^{-1}(z)$, it gives  \begin{eqnarray*}R_\mu=C_{\sqrt{\mu\bxt\mu_\la}}&\iff&\lf(z\times T\circ \lf(\f{M_\mu^{-1}(z)}{1+\la z}\ri)^{-1}\ri)^{-1}\circ ((z+1)M_\mu^{-1}(z))\times  T(z) =(z+1)M_\mu^{-1}(z).\end{eqnarray*}Dividing by $T(z)$, it gives \begin{eqnarray*}R_\mu=C_{\sqrt{\mu\bxt\mu_\la}}&\iff&\lf(z\times T\circ \lf(\f{M_\mu^{-1}(z)}{1+\la z}\ri)^{-1}\ri)^{-1}\circ ((z+1)M_\mu^{-1}(z)) =\f{M_\mu^{-1}(z)}{1+\la z}\\ &\iff& (z+1)M_\mu^{-1}(z) =\lf(z\times T\circ \lf(\f{M_\mu^{-1}(z)}{1+\la z}\ri)^{-1}\ri)\circ\f{M_\mu^{-1}(z)}{1+\la z}\\ &\iff& (z+1)M_\mu^{-1}(z) = \f{M_\mu^{-1}(z)}{1+\la z}T(z),
\end{eqnarray*}which is obviously true by definition of $T(z)$.
\end{pr}

\subsection{Main result of the paper}

The main theorem of this paper is the following one. $\la\in [0,1]$ is still fixed.
\begin{Th}\label{27.08.08.1}For any pair $\mu,\mu'$ of \pro measures on $[0,+\infty)$, we have 
\be\label{26.08.08.2}\sqrt{\mu\bxt\mu_\la}\arc \sqrt{\mu'\bxt\mu_\la}=\sqrt{(\mu\bxp\mu')\bxt\mu_\la}.\ee\end{Th}

\begin{rmq}\label{28.08.08.1avant} {\rm Note that in the case where $\la=0$,  this theorem expresses what we already knew about $\bxp_0$ (and which is explained in  the second paragraph of the introduction), but that the case $\la=1$ isn't a consequence of the already known formula $\bxp_1=\bxp$.}
 \end{rmq}

\begin{rmq}\label{28.08.08.1}{\rm Part of this theorem could have been deduced from Theorem 6 of \cite{dr07}. However, \eqref{26.08.08.2} could be deduced from the theorem of Debbah and Ryan only for laws $\mu,\mu'$ which can be expressed as limit singular laws of $n$ by $p$ (for $n/p\simeq \la$) corners of large $p\times p$ bi-unitarily invariant random matrices, but it follows from Theorem 14.10 of \cite{ns06} that not every law has this form. Moreover, even though the idea which led us to our result was picked in the  pioneer work of Debbah and Ryan, our proof is much shorter and shows the connection with the rectangular machinery in a more clear way (via Theorem \ref{R=C} and Remark \ref{28.08.08.16}).} \end{rmq}

\begin{pr}Define $\nu:=\sqrt{\mu\bxt\mu_\la}\arc \sqrt{\mu'\bxt\mu_\la}$. By \eqref{28.08.08.17}, we have  $$C_\nu=C_{\sqrt{\mu\bxt\mu_\la}}+ C_{\sqrt{\mu'\bxt\mu_\la}}.$$ Thus, by Theorem \ref{R=C}, and \eqref{25.08.08.1}, we have $$C_\nu=R_\mu+R_{\mu'}=R_{\mu\bxp\mu'}=C_{\sqrt{(\mu\bxp\mu')\bxt\mu_\la}}.$$Hence by injectivity of the rectangular $R$-transform (Theorem 3.8 of \cite{bg07}), \eqref{26.08.08.2} is valid.\end{pr}

The formula \eqref{26.08.08.2} gives us a new insight on rectangular free convolutions: it allows to express it, in certain cases, in terms of the free convolutions ``of  square type" $\bxp$ and $\bxt$. However, only laws which can be expressed under the form\be\label{27.08.08.6}\sqrt{\mu\bxt\mu_\la},\quad\textrm{ ($\mu$ \pro measure on $[0,+\infty)$)}\ee  can have their rectangular convolution computed via formula \eqref{26.08.08.2}. Thus it seems natural to ask whether all symmetric laws %try to find out which laws 
can be expressed  like in \eqref{27.08.08.6}. Note that it is equivalent to the fact that any law on $[0,+\infty)$ can be expressed under the form $\mu\bxt \mu_\la$, which is equivalent to the fact that the Dirac mass at one $\delta_1$  can be expressed under the form $\mu\bxt \mu_\la$. Indeed, if $\delta_1=\mu\bxt\mu_\la$, then  any law $\tau$ on $[0,+\infty)$ satisfies $\tau=\tau\bxt\delta_1=(\tau\bxt\mu)\bxt\mu_\la$. The following proposition shows that it is not the case. However, Theorem \ref{27.08.08.5} will show that many symmetric laws can be expressed like in \eqref{27.08.08.6}.

\begin{propo}\label{27.08.08.9}Unless $\la=0$, the law  $\f{\delta_1+\delta_{-1}}{2}$ cannot be expressed  under the form $\sqrt{\mu\bxt\mu_\la}$ for $\mu$ \pro measure on $[0,+\infty)$. \end{propo}
 
 \begin{pr} Suppose that $\la>0$ and that there is a \pro measure $\mu$ on $[0,+\infty)$ \st $\f{\delta_1+\delta_{-1}}{2}=\sqrt{\mu\bxt\mu_\la}$. Then $\delta_1=\mu\bxt\mu_\la$. This is impossible,  by Corollary  3.4 of \cite{b06}, which states that the free multiplicative convolution of two laws which are not Dirac masses has always a non null absolutely continuous part (there is another, more direct way to see that it is impossible: by \eqref{25.08.08.1}, such a law $\mu$ has to satisfy $S_\mu(z)=1+\la z$, which implies that for $z$ small enough, $M_\mu(z)=\f{z-1+[(1-z)^2+4\la z]^{1/2}}{2\la}$: such a function doesn't admit any analytic continuation to $\C\bck[0,+\infty)$, thus no such \pro measure $\mu$ exists).
  \end{pr}

Theorem \ref{27.08.08.1} has a consequence on the free convolutions ``of  square type" which wasn't known yet, despite the many papers written the last years about questions related to the arithmetics of  these convolutions, e.g.  \cite{bv95,appenice,fbg04,gotzea, gotzeb,bbg07,bbcc08}.
\begin{cor}\label{27.08.08.1cor}For any pair $\mu,\mu'$, of \pro measures on $[0,+\infty)$ we have 
\be\label{27.08.08.22.10.6.09}
\sqrt{\mu\bxt\mu_1}\bxp \sqrt{\mu'\bxt\mu_1}=\sqrt{(\mu\bxp\mu')\bxt\mu_1}.\ee \end{cor}

\begin{pr}It is an obvious consequence of Theorem \ref{27.08.08.1} and of the fact that $\bxp_1=\bxp$.
\end{pr}

\begin{rmq}\label{10.6.09.1}{\rm The referee of the paper communicated to us a proof of \eqref{27.08.08.22.10.6.09} which is not, as ours, based on computations on the $R$- and $S$-transforms, but   on the direct proof of \eqref{10.6.09.2} in the special case $\la=1$. Let us briefly outline  this proof.  When $\la=1$, by \cite[Eq. (4.1)]{fbg05.inf.div}, \eqref{10.6.09.2} reduces to \be\label{bayrou-cohn.10.6.09}k_n(\mu)=k_{2n}(\sqrt{\mu\bxt\mu_1}).\ee Let  $a,s$ are free elements in a tracial non commutative \pro space %$(\mc{A},\vfi)$, 
with respective distributions $\mu$ and the standard semicircle law. By  \cite[Prop. 12.13]{ns06}, $s^2$ has distribution $\mu_1$,  hence $sas$ has distribution $\mu\bxt\mu_1$. It follows, by \cite[Prop. 12.18]{ns06}, that for all $n$, the $n$-th moment of $\mu$ is equal to $k_n(\mu\bxt\mu_1).$ But by \cite[Prop. 11.25]{ns06}, for all $n$, we have $$k_n(\mu\bxt\mu_1)=\sum_{\pi\in \NC(n)}\prod_{V\in \pi} k_{2|V|}(\sqrt{\mu\bxt\mu_1}).$$ It follows, using the expression of the $n$-th moment of $\mu$ in terms of its free cumulants, that for all $n$, $$\sum_{\pi\in \NC(n)}\prod_{V\in \pi} k_{|V|}(\mu)=\sum_{\pi\in \NC(n)}\prod_{V\in \pi} k_{2|V|}(\sqrt{\mu\bxt\mu_1}),$$ and that for all $n$, $k_n(\mu)=k_{2n}(\sqrt{\mu\bxt\mu_1})$.}\end{rmq}

\section{Consequences on square and rectangular infinite divisibility}
\subsection{Prerequisites on infinite divisibility and L\'evy-Kinchine formulas}
 Infinite divisibility is a fundamental probabilistic notion, at the base of L\'evy processes, and which allows to explain deep relations between limit theorems for sums of either independent random variables, square or rectangular random matrices. Let us briefly recall basics of this theory \cite{gne,sato,defconv, appenice,fbg05.inf.div}.

Let $*$ denote the classical convolution of \pro measures on the real line.
 Firstly, recall that a \pro measure $\mu$ is said to be $*$-infinitely divisible (resp.  $\bxp$-, $\arc$-infinitely divisible) if for all integer $n$, there exists a \pro measure $\nu_n$ \st $\nu_n^{*n}=\mu$ (resp.   $\nu_n^{\bxp n}=\mu$,  $\nu_n^{\arc n}=\mu$). In this case, there exists a $*$- (resp. $\bxp$-, $\arc$-) semigroup $(\mu_t)_{t\geq 0}$ \st $\mu_0=\delta_0$ and $\mu_1=\mu$. For all $t$, $\mu_t$ is denoted by $\mu^{* t}$ (resp. $\mu^{\bxp t}, \mu^{\arc t}$). Infinitely divisible  distributions have been classified: $\mu$ is $*$- (resp. $\bxp$-) infinitely divisible \ssi there exists a real number $\gamma$ and a positive finite measure on the real line   $\sigma$ \st  the Fourier transform is $\hat{\mu}(t)=\exp\lf[i\gamma t+\int_\R(e^{itx}-1-\f{itx}{x^2+1})\f{x^2+1}{x^2}\ud \sigma(x)\ri]$ (resp. $R_\mu(z)=\gamma z+z\int_\R\f{z+t}{1-tz}\ud \sigma(t)$). Moreover, in this case, such a pair $(\gamma, \sigma)$ is unique, it is called the {\it L\'evy pair} of $\mu$ and we denote $\mu$ by $\nu_{*}^{\gamma, \sigma}$ (resp. $\nu_{\bxp}^{\gamma, \sigma}$). For all $t\geq 0$,  $\mu_t$ has L\'evy pair $(t\gamma, t\sigma)$. In the same way, a symmetric \pro measure $\nu$ is  $\arc$-infinitely divisible \ssi there exists a positive finite symmetric measure on the real line $G$ \st $C_\nu(z)=z\int_{\R}\f{1+t^2}{1-zt^2}\ud G(t)$. In this case, the measure $G$ is unique, and $\nu$ will be denoted by $\nu_{\arc}^G$. The correspondences $\nu_{*}^{\gamma, \sigma}\longleftrightarrow \nu_{\bxp}^{\gamma, \sigma}$ (for any pair $(\gamma, \sigma)$ as above) and $\nu_{*}^{0, G}\longleftrightarrow \nu_{\arc}^{G}$ (for any $G$ as above) are called {\it Bercovici-Pata bijections}. These bijections have many deep properties \cite{appenice, fbg05.inf.div}, some of which will be mentioned in the proof of the following lemma.

\begin{lem}\label{inf.div.pos.08.08}
Let $\gamma$ be  a real number and $\sigma$ be a positive finite measure on the real line. Then we have equivalence between: 

(i) For all $t\geq 0$, $\nu_{*}^{t\gamma, t\sigma}$ is supported on $[0,+\infty)$.

(ii) For all $t\geq 0$, $\nu_{\bxp}^{t\gamma, t\sigma}$ is supported on $[0,+\infty)$.

(iii) We have $\sigma((-\infty,0])=0$ and  the integral $\int \ff{x}\ud \sigma(x)$ is finite and $\leq \gamma$.
\end{lem}

\begin{pr}The equivalence between (i) and (iii) follows from Theorem 24.7 and Corollary 24.8 of \cite{sato}. Let us prove the equivalence between (i) and (ii). In order to do that, let us recall a fact proved in \cite{appenice}: for any L\'evy pair $(\gamma, \sigma)$ and any sequence $(\nu_n)$ of \pro measures, one has \be\label{27.08.08.3}\nu_n^{*n}\textrm{ converges weakly to }\nu_{*}^{\gamma, \sigma}\iff\nu_n^{\bxp n}\textrm{ converges weakly to }\nu_{\bxp}^{\gamma, \sigma}.\ee Let us suppose (i) (resp. (ii)) to be true. Let us fix $t\geq 0$. For all $n$, we have $$(\nu_{*}^{\f{t\gamma}{n}, \f{t\sigma}{n}})^{*n}=\nu_{*}^{t\gamma, t\sigma}\quad\textrm{ (resp. }(\nu_{\bxp}^{\f{t\gamma}{n}, \f{t\sigma}{n}})^{\bxp n}=\nu_{\bxp}^{t\gamma, t\sigma}\textrm{).}$$Thus by \eqref{27.08.08.3}, $$(\nu_{*}^{\f{t\gamma}{n}, \f{t\sigma}{n}})^{\bxp n}\textrm{ converges weakly to }\nu_{\bxp}^{t\gamma, t\sigma}\quad\textrm{ (resp. }(\nu_{\bxp}^{\f{t\gamma}{n}, \f{t\sigma}{n}})^{* n}\textrm{ converges weakly to }\nu_{*}^{t\gamma, t\sigma}\textrm{).}$$Thus since any free (resp. classical) additive convolution and any weak limit of measures with supports on $[0,+\infty)$ has support on $[0,+\infty)$, (ii) (resp. (i)) holds.
\end{pr}

\begin{rmq}\label{27.08.08.55}{\rm Note that (i) is equivalent to the fact that there exists $t> 0$ \st $\nu_*^{t\gamma, t\sigma}$ is supported on $[0,+\infty)$  \cite[Cor. 24.8]{sato}. However,  the same is not true for the free infinitely divisible laws. Indeed, let, for each $t\geq 0$, $\MP_t$ denote the   Marchenko-Pastur law with mean $t$ \cite[Ex. 3.3.5]{hiai} and let us define, for each $t$, $\mu_t=\MP_t*\delta_{-t/4}$. Then since free and classical convolutions with Dirac masses are the same, $(\mu_t)_{t\geq 0}$ is a convolution semi-group with respect to $\bxp$. But  $\mu_4$ is supported on $[0,+\infty)$, whereas for each $t\in(0,1]$, the support of $\mu_t$ contains a negative number (namely $-t/4$).}
\end{rmq}

\subsection{Main result of the section}The following theorem allows us to claim that even though not every symmetric law can be expressed under the form $\sqrt{\mu\bxt\mu_\la}$ for $\mu$ law on $[0,+\infty)$ (see Proposition \ref{27.08.08.9}), many of them have this form. $\la\in [0,1]$ is still fixed. 

For $G$ measure on the real line, we let $G^2$ denote the push-forward of $G$ by the function $t\mapsto t^2$.

\begin{Th}\label{27.08.08.5}(i) Let $\mu$ be a $\bxp$-infinitely divisible law \st for all $t\geq 0$, $\mu^{\bxp t}$ is supported on $[0,+\infty)$. Then the law $\sqrt{\mu\bxt\mu_\la}$ is $\arc$-infinitely divisible, with L\'evy measure the only symmetric measure $G$ \st \be\label{27.08.08.33}G^2=\lf(\gamma-\int\f{1}{x}\ud \sigma(x)\ri)\delta_0+ \f{1+x^2}{x(1+x)}\ud \sigma(x), \ee where  $(\gamma, \sigma)$ denotes the  L\'evy pair  of $\mu$.

(ii) Reciprocally, any $\arc$-infinitely divisible law $\nu$ has the form $\sqrt{\mu\bxt\mu_\la}$ for some $\bxp$-infinitely divisible law $\mu$ \st for all $t\geq 0$, $\mu^{\bxp t}$ is supported on $[0,+\infty)$. Moreover, the L\'evy pair  $(\gamma, \sigma)$ of $\mu$ is defined by \be\label{27.08.08.44}\gamma=\int_{[0,+\infty)}\f{1+x}{1+x^2}\ud G^2(x)\quad\textrm{ and }\quad \sigma=\f{x(1+x)}{1+x^2}\ud G^2(x),\ee where $G$ denotes the L\'evy measure of $\nu$.
\end{Th}

\begin{pr} (i) Note that by Theorem \ref{27.08.08.1}, the map $\mu\mapsto\sqrt{\mu\bxt\mu_\la}$ is a morphism from the set of laws on $[0,+\infty)$ to the set on symmetric laws on the real line endowed respectively with the operations $\bxp$ and $\arc$. Thus if $\mu$ is $\bxp$-infinitely divisible, then $\sqrt{\mu\bxt\mu_\la}$ is $\arc$-infinitely divisible. Moreover, if the L\'evy pair of $\mu$ is $(\gamma, \sigma)$, then its $R$-transform is $R_\mu(z)=\gamma z+z \int_{t\in \R} \f{z+t}{1-zt}\ud \sigma(t).$ By Theorem \ref{R=C}, it implies that $C_{\sqrt{\mu\bxt\mu_\la}}(z)= \gamma z+z \int_{t\in \R} \f{z+t}{1-zt}\ud \sigma(t)$. 
But by uniqueness, the L\'evy measure $G$ of $\sqrt{\mu\bxt\mu_\la}$ is characterized by the fact that 
  $C_{\sqrt{\mu\bxt\mu_\la}}(z)=z\int_\R\f{1+t^2}{1-zt^2}\ud G(t)$.   Thus to prove \eqref{27.08.08.33}, it suffices to prove that for $G$ given by  \eqref{27.08.08.33}, for all $z$, one has $$\gamma z+z \int_{t\in \R} \f{z+t}{1-zt}\ud \sigma(t)=z\int_\R\f{1+t^2}{1-zt^2}\ud G(t),$$which can easily be verified.

(ii) Let $\nu$ be a $\arc$-infinitely divisible law with L\'evy measure denoted by $G$. Let $(\gamma, \sigma)$ be the L\'evy pair defined by \eqref{27.08.08.44}. Note that $(\gamma, \sigma)$ satisfies (iii) of Lemma \ref{inf.div.pos.08.08}, thus, for $\mu:=\nu_{\bxp}^{\gamma,\sigma}$, for  all $t\geq 0$,  the law $\mu^{\bxp t}$ is actually supported by $[0,+\infty)$. Thus by (i), $\sqrt{\mu\bxt\mu_\la}$ is $\arc$-infinitely divisible with L\'evy measure the only symmetric measure $H$ satisfying $$H^2=\lf(\gamma-\int\f{1}{x}\ud \sigma(x)\ri)\delta_0+ \f{1+x^2}{x(1+x)}\ud \sigma(x).$$ To prove that  $\sqrt{\mu\bxt\mu_\la}=\nu$, it suffices to prove that $H=G$, which can easily be verified.
\end{pr}

One of the consequences of this theorem is that it gives us a description of the free multiplicative  convolution of two Marchenko-Pastur laws (i.e. free Poisson laws), one of them having a mean $\geq 1$. For all $t>0$, the Marchenko-Pastur law  $\MP_t$ with mean $t$ has been introduced  at Remark \ref{27.08.08.55}.

\begin{cor}\label{28.08.08.randocycle}Consider $a,c>0$ \st $a>1$. Then $\MP_c\bxt\MP_a$ is the push forward, by the map $x\mapsto ax^2$, of the $\arc$-infinitely divisible law with L\'evy measure $\f{c}{4}(\delta_1+\delta_{-1})$ for $\la=1/a$. 
\end{cor}
 
 \begin{pr}
 It suffices to notice that for $\la=1/a$, $\MP_a$ is the push-forward, by the map $x\mapsto ax$, of the law $\mu_\la$, that $\MP_c$ is the $\bxp$-infinitely divisible law with L\'evy pair $(c/2,c/2\delta_1)$, and then to apply (i) of Theorem \ref{27.08.08.5}.
 \end{pr}
 
\begin{figure}[h!]
\begin{center}
%\rotatebox{270}{\scalebox{0.6}{
\scalebox{0.6}{\includegraphics{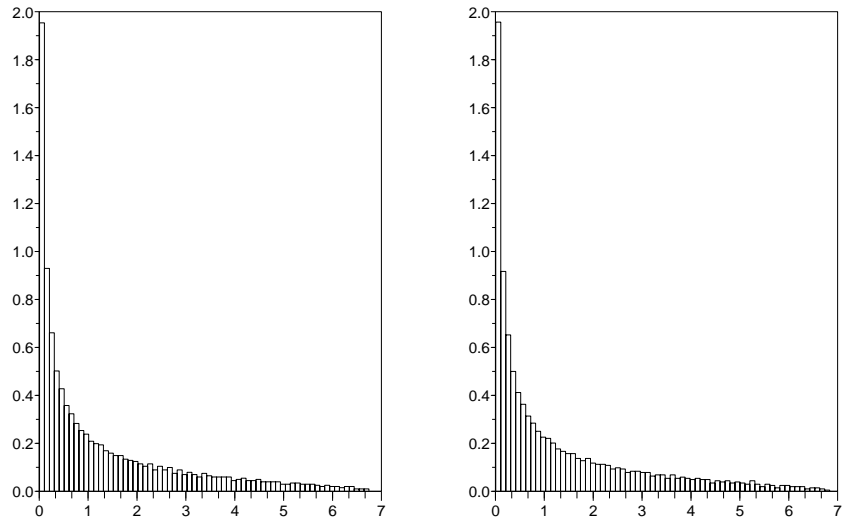}}
\caption{Histograms  of the spectrums of $MM^*$ (left) and $\ff{nq} TT^*QQ^*$ (right) for $n=2000$, $\la=0.6$, $c=1.3$.}\label{histo.11.06.09}
\end{center}
\end{figure}

 This corollary can be interpreted as the coincidence of the limit laws of two different matrix models. Indeed, the  $\arc$-infinitely divisible law with L\'evy measure $\f{c}{4}(\delta_1+\delta_{-1})$ was already known \cite[Prop. 6.1]{fbg05.inf.div} to be the limit symmetrized singular law of the random matrix $M:=\sum_{k=1}^pu_kv_k^*$, for  $n,p,q$ tending to infinity in such a way that $p/n\to c$  and $n/q\to \la$ and  $(u_k)_{k\geq 1}$, $(v_k)_{k\geq 1}$ two independent families  of independent random vectors \st for all $k$, $u_k, v_k$ are  uniformly distributed on the unit spheres of respectively $\C^n, \C^q$. 
Thus, if, for large $n,p,q$'s \st $p/n\simeq c$ and $n/q\simeq \la$, one considers such a random matrix $M$ and also
  two independent random matrices $T,Q$  with respective dimensions $n \times p, n\times q$, the   entries of which   are  independent real standard Gaussian random variables, then the empirical spectral measures of 
  the random matrices $MM^*$ and $\ff{nq} TT^*QQ^*$ are close to each other, as illustrated by Figure \ref{histo.11.06.09}. 
    
\end{document}